\documentclass[12pt]{amsart}
\theoremstyle{plain}
\newtheorem{Thm}{Theorem}

\newtheorem{Prop}[Thm]{Proposition}

\newtheorem{Def}[Thm]{Definition}

\begin{document}

\title[Gradient Ricci soliton]
{Remarks on Gradient Ricci Solitons}

\author{Li MA }

\address{Department of mathematical sciences \\
Tsinghua university \\
Beijing 100084 \\
China}

\email{lma@math.tsinghua.edu.cn} \dedicatory{}
\date{Oct 20th, 2004}

\thanks{$^*$ This work is supported in part by
the Key 973 project of Ministry of Science and Technology of
China}

\begin{abstract}
In this paper, we study the gradient Ricci soliton equation on a
complete Riemannian manifold. We show that under a natural decay
condition on the Ricci curvature, the Ricci soliton is Ricci-flat
and ALE.
\end{abstract}

 \maketitle

\section{Introduction}
In the study of Ricci flow on a Riemnnian manifold, we meet Ricci
solitons (see \cite{H95}). Ricci solitons and their sisters
Kahler-Ricci solitons are important objects by their own right
(see also \cite{B04}, \cite{C96}, and \cite{T04} ).

  Given a Riemannian manifold $(X,g)$ of dimension $n$. Let $Rc$
  be the Ricci tensor of the metric $g$. The equation for a homothetic Ricci
soliton is
$$
Rc=cg+{L}_Vg
$$
where $c$ is a homothetic constant, $V$ is a smooth vector field
on $X$, and ${L}_Vg$ is the Lie derivative of the metric $g$. When
$c=0$, the soliton is steady. For $c>0$ the soliton is shrinking,
and one can consider the Ricci flow on the sphere as such an
example. For $c<0$ the soliton is expanding. When $V$ is the
gradient of a smooth function, we call such solitons  {\em
Gradient Homothetic Ricci Solitons}. Let $R$ be the scalar
curvature of $g$. An important equation in Riemannian geometry and
general relativity theory is the so called Einstein equation:
$$
E_{ij}=T_{ij},
$$
where $E=Rc-\frac{R}{2}g$, and $T$ is the energy momentum tensor
in the space. The tensor $T$ sometimes is also an unknown being.
So it is interesting to know whether $T$ is the Hessian matrix of
a smooth function, and to explore more properties about this
function.

  In this short paper we constrict ourselves to consider the following
 question about the gradient Ricci soliton on $X$.
  We deduce a Liouville type theorem for a smooth
  solution $f$ for the equation
  $$
Rc=D^2f
  $$
  on $X$,
  where $D^2f$ is the Hessian matrix of the function $f$.

  Taking
  the trace of the equation above, we get that
$$
R=\Delta f.
$$
Here we denote by $R$ the scalar curvature of the metric
  $g$ and $\Delta f$ the Laplacian of $f$. If such a smooth function exists,
  we call $(X,g)$ a {\em gradient Ricci soliton}.

  In the
  following, we write $D$ the covariant derivative of $g$.

  As in \cite{H95}, we can derive (see the next section) that there is a constant $M$
  such that
  $$
|Df|^2+R=M
  $$
  on $X$. Then we have the following simple observation.

  \begin{Prop} Assume $X$ is compact. Then $(X,g)$ is Ricci-flat.
  \end{Prop}
\begin{proof}
 In fact, we have
 $$
  \Delta f+|Df|^2=M.
 $$
 Set
 $$
u=e^f.
 $$
 Then we have that $u>0$ and
 $$
\Delta u=Mu
 $$
 Integrating we find
$$ M\int_X u=0.
$$
Then we have $M=0$. Hence, $u$ is a harmonic function. By the
Maximum principle, we know that $u$ is a constant, so $f$ is a
constant. This implies that $Rc=0$ and $(X,g)$ is Ricci-flat.
\end{proof}

  We want to generalize the above result. Observe that if $(X,g)$
  is non-compact and complete, the matter is not simple.
  However, if $R(x)\to 0$ as $r:=r(x)\to\infty$, then we can get some conclusions.
  Here
  $r(x)=dist(x,o)$ is the distance of the point $x$ from a fixed point $o$.

  Since $M\geq R$ and $R$ can be very small at infinity, we have
  $M\geq 0$. If we also assume that $|Df(x)|\to 0$ as $r(x)\to
  \infty$, then we must have that $M=0$. Hence $u$ is again a positive
  harmonic function.
   If we further assume that $Rc\geq
  0$, then we must have by the Liouvile theorem of Yau \cite{Y} that
$u$ is a constant.
 So we have obtained the following result.

\begin{Prop} Assume $X$ is a non-compact complete Riammnian manifold
with non-negative Ricci curvature. If either case 1) $\int_X
|R|^p<+\infty$ and $\int_X |f|^p<+\infty$ or $\int_X
|Df|^p<+\infty$ ( for some $p\geq 1$) or case 2) both $R$ and $f$
decay to zero at infinity, then $f$ is a constant and $(X,g)$ is
Ricci-flat.
  \end{Prop}
  \begin{proof} In both cases we can easily conclude that $M=0$.
 Then we have that
$$
\Delta u=M=0.
$$
By the Liouvile theorem of Yau \cite{Y} we have that $u$ is a
constant.

  \end{proof}

Now the nontrivial matter is to treat the case when $Rc$ has no
sign assumption. The natural consideration is sing the cut-off
function trick. In fact, it works for some cases. We have the
following theorem, which is the main result of this paper.

\begin{Thm}
Assume $X$ is a non-compact complete Riammnian manifold which is
quasi-isometric to Euclidean space at infinity. If $n\geq 3$ and
$Rc$ satisfies that $$\int_X|Rc|^{n/2}<+\infty ,$$ then $f$ is a
constant and $(X,g)$ is Ricci-flat and ALE of order $n-1$. If we
further assume that $n=4$, then $(X,g)$ is ALE of order $4$.
\end{Thm}

The definition about ALE space will be given in the next section.
In our proof of the result above, we will use the Bochner formula,
Moser iteration method, and  the interesting result of
Bando-Kasue-Nakajima \cite{BKN89}. Our idea can also be used to
study Ricci-Kahler soliton and related question for Einstein
equations.

\section{Notations, definitions, and basic facts}

In local coordinates $(x^i)$ of the Riemannian manifold $(X,g)$,
we write the metric $g$ as $(g_{ij})$. The corresponding
Riemannian curvature tensor and Ricci tensor are denoted by
$Rm=(R_{ijkl})$ and $Rc=(R_{ij})$ respectively. Hence,
$$
R_{ij}=g^{kl}R_{ikjl}
$$
and
$$
R=g^{ij}R_{ij}.
$$
We write the covariant derivative of a smooth function $f$ by
$Df=(f_i)$, and denote  the Hessian matrix of the function $f$ by
$D^2f=(f_{ij})$, where  $D$ the covariant derivative of $g$ on
$X$.
 The
higher order covariant derivatives are denoted by $f_{ijk}$, etc.
Similarly, we use the $T_{ij,k}$ to denote the covariant
derivative of the tensor $(T_{ij})$. We write
$T_j^i=g^{ik}T_{jk}$. Then the Ricci soliton equation is
$$
R_{ij}=f_{ij}.
$$
Taking covariant derivative we get
$$
f_{ijk}=R_{ij,k}.
$$
So we have
$$
f_{ijk}-f_{ikj}=R_{ij,k}-R_{ik,j}.
$$
By the Ricci formula we have that
$$
f_{ijk}-f_{ikj}=R_{ijk}^lf_l.
$$
Hence we obtain that
$$
R_{ij,k}-R_{ik,j}=R_{ijk}^lf_l.
$$
Recall that the contracted Bianchi identity is
$$
R_{ij,j}=\frac{1}{2}R_i.
$$
Upon taking the trace of the previous equation we get that
$$
\frac{1}{2}R_i+R_{i}^kf_k=0,
$$
i.e.,
$$
R_k=-2R_{k}^jf_j.
$$
Then
$$
D_k(|Df|^2+R)=2f_j(f_{jk}-R_{jk})=0.
$$
So, $|Df|^2+R$ is a constant, which is denoted by $M$.

 To analyze the global behavior of the geometry of the manifold
 $X$,
 we need the following definition of quasi-isometry.

 \begin{Def}
We say that $(X,g)$ is quasi-isometric to the Euclidean space
$R^n$ if there is a compact subset $K$ in $X$ such that the metric
$X-K$ is uniformly equivalent to $R^n-B_1(0)$.
 \end{Def}

It is well-known that for such a Riemannian manifold, the Sobolev
inequality is true. We shall use this fact to make the Moser
iteration method work in $X$. A special case of such a manifold
with one end  is called {\em asymptotically locally Euclidean}
(ALE) of order $\tau>0$ (see \cite{B86}); namely there are
constant $\rho>0$, $\alpha \in (0,1)$ and a
$C^{\infty}$-diffeomorphism
$$\Xi: X-K\to R^n-B_{\rho}(0)$$ such that $\phi=\Xi^{-1}$ satisfies
$$
(\phi^{*}g)_{ij}(x)=\delta_{ij}+ 0(|x|^{-\tau}),
$$
$$
\partial_k(\phi^{*}g)_{ij}(x)=0(|x|^{-1-\tau}),
$$
and
$$
\frac{|\partial_k(\phi^{*}g)_{ij}(x)-\partial_k(\phi^{*}g)_{ij}(z)|}{|x-z|^{\alpha}}
=0(|x|^{-1-\alpha-\tau},|z|^{-1-\alpha-\tau})
$$
for any $x,z\in R^n-B_{\rho}(0)$.

\section{Proof of Theorem 3}
  We now recall the Bochner formula for any smooth function $w$ on $X$:
  $$
\frac{1}{2}\Delta|D w|^2=|D^2w|^2+<D\Delta w,Dw>+Rc(Dw,Dw).
  $$
For our function $f$ in the Ricci soliton equation, we have that
$$
\frac{1}{2}\Delta|Df|^2=|D^2f|^2+<DR,Df>+Rc(Df,Df)
$$
Recall that
$$
R_k=-2R_{kj}f_j.
$$
Then we have
$$ <DR,Df>=-2Rc(Df,Df).
$$
 So we have
$$
\frac{1}{2}\Delta|Df|^2=|Rc|^2-Rc(Df,Df).
$$
This equation can also be written as
$$
-\frac{1}{2}\Delta R=|Rc|^2-\frac{1}{2}(DR,Df),
$$
 from which we know the following inequality
$$
\frac{1}{2}\Delta|Df|^2\geq-|Rc||Df|^2.
$$

 We shall use this differential inequality to study the
behavior of $|Df|$ at infinity. Let $w=|Df|^2$. Then $w=|R|$.
Using the assumption that
$$ \int
|Rc|^{n/2}<+\infty
$$
we get that
$$ \int
|w|^{n/2}<+\infty.
$$
Using Moser's iteration method
 we know (see Lemmas 4.1,4.2,4.6 in \cite{BKN89}) that there is a
 positive constant $\alpha>0$ such that
 $$
w=0(r^{-\alpha}).
 $$
That is to say that
$$
|R(x)|=|Df|^2(x)=0(r^{-\alpha}).
$$
This implies that $M=0$, and $R\leq 0$ on $X$. Using the maximum
principle to the equation
\begin{equation}
-\frac{1}{2}\Delta R=|Rc|^2-\frac{1}{2}(DR,Df)
\end{equation}
and the fact that $$ |Rc|^2\geq \frac{1}{n}R^2,
$$
we find that the minimum of $R$ on $X$ can not be negative, and it
must be zero. Hence $R=0$ on $X$. This implies from (1) that
$Rc=0$. Using Theorem 1.5 in \cite{BKN89} we know that $(X,g)$ is
ALE of order $n-1$. If further $n=4$, $(X,g)$ is ALE of order $4$.

This proves our Theorem 3.

{\bf Acknowledgement}. The author thanks Mr. D.Z. Chen for
checking spelling mistakes in a previous version of this paper.


\begin{thebibliography}{20}
\bibitem{BKN89}
S.Bando, A.Kasue, and H.Nakajima, {\em On a construction of
coordinates at infinity on manifolds with fast curvature decay and
maximal volume growth}, Invent. Math., 97(1989)313-349.

\bibitem{B86}
Bartnik, {\em The mass of an asymptotically flat manifold}, Comm.
Pure Appl. Math., 34(1986)661-693.

\bibitem{B04}
R.Bryant, {\em Gradient Kahler-Ricci solitons},
ArXiv.math.DG/0407453

\bibitem{C96}
H.D.Cao, {\em Existence of gradient Kahler-Ricci soliton},
Elliptic and parabolic methods in geometry, Eds. B.Chow,
R.Gulliver, S.Levy, J.Sullivan, A K Peters, pp.1-6, 1996.

\bibitem{H95}
 R.Hamilton, {\em The formation of Singularities in the Ricci flow},
 Surveys in Diff. Geom.,
 Vol.2, pp7-136, 1995.

\bibitem{T04}
G.Tian, Lectures on some geometric problems, Tsinghua University,
2004.

\bibitem{Y}
S.T.Yau, {\em Harmonic functions on complete Riemannian
manifolds}, Comm. Pure Appl. Math., 28(1975)201-228.

\end{thebibliography}
\end{document}